\documentclass{amsart}

\usepackage{amsmath,amssymb,amsthm,pinlabel}

\def\co{\colon\thinspace}
\newcommand{\Int}{\mbox{\rm Int}}
\newcommand{\N}{\mathbb{N}}
\newcommand{\R}{\mathbb{R}}
\newcommand{\Z}{\mathbb{Z}}
\newcommand{\LL}{\mathbb{L}}
\newcommand{\oL}{\overline{L}}
\newcommand{\tb}{{\tt tb}}
\newcommand{\otb}{\overline{\tt tb}}
\newcommand{\lk}{{\tt lk}}
\newcommand{\rot}{{\tt rot}}

\newtheorem{thm}{Theorem}

\newtheorem{prop}[thm]{Proposition}
\newtheorem{cor}[thm]{Corollary}

\theoremstyle{definition}
\newtheorem*{rem}{Remark}

\newtheorem*{ack}{Acknowledgements}

\begin{document}

\title{Legendrian helix and cable links}

\author{Fan Ding}
\address{Department of Mathematics, Peking University,
Beijing 100871, P.~R. China}
\email{dingfan@math.pku.edu.cn}
\author{Hansj\"org Geiges}
\address{Mathematisches Institut, Universit\"at zu K\"oln,
Weyertal 86--90, 50931 K\"oln, Germany}
\email{geiges@math.uni-koeln.de}
\date{}

\begin{abstract}
Lisa Traynor has described an example of a two-component Legendrian
`circular helix link' $\Lambda_0\sqcup\Lambda_1$ in the $1$--jet space
$J^1(S^1)$ of
the circle (with its canonical contact structure) that is topologically
but not Legendrian isotopic to the link $\Lambda_1\sqcup\Lambda_0$. We give
a complete classification of the Legendrian realisations of this
topological link type, as well as all other `cable links' in~$J^1(S^1)$.
\end{abstract}

\maketitle

\section{Introduction}
Considerable progress has been made towards the classification
of Legendrian knots and links in contact $3$--manifolds, e.g.\
\cite{etho01,etho03,ngtr04,dige07}. For a general
introduction to this topic see~\cite{etny05}. In~\cite{tray97},
Lisa Traynor exhibited an intriguing example
of an ordered Legendrian two-component link that is topologically
but not Legendrian isotopic to the link obtained by interchanging
the two components. The purpose of the present paper
is to give a complete classification of the Legendrian realisations
of this and some related topological link types.

Let
\[ J^1(S^1)=T^*S^1\times\R=\{ (x,y,z)\co x\in\R /2\pi\Z ,\, y,z\in\R\} \]
be the $1$--jet space of $S^1$ with its standard contact structure
\[ \xi =\ker (dz-y\, dx ).\]
Here $y$ denotes the fibre coordinate in $T^*S^1$, and $z$ the
coordinate in the $\R$--factor of $J^1(S^1)$. The graph of any smooth
function $g\co S^1\rightarrow\R$ lifts to a Legendrian knot
\[ \Lambda_g:=\{ (x,y,z)\co y=g'(x),\, z=g(x)\} .\]
The linear interpolation between two smooth functions $g,h\co S^1\rightarrow
\R$ gives rise to a Legendrian isotopy between $\Lambda_g$ and
$\Lambda_h$. The following theorem, therefore, is rather surprising.
Here we write $\Lambda_s$ for the Legendrian knot corresponding to
the constant function of value $s\in\R$, oriented by the variable~$x$.

\begin{thm}[Traynor]
\label{thm:Traynor}
The ordered Legendrian links $\Lambda_0\sqcup\Lambda_1$ and
$\Lambda_1\sqcup\Lambda_0$ in $(J^1(S^1),\xi )$ are topologically, but
not Legendrian isotopic.
\end{thm}

In Theorem~\ref{thm:helix} below we give a complete classification
of the Legendrian realisations of the `circular helix link',
i.e.\ the topological link type
$\Lambda_0\sqcup\Lambda_1$ in $(J^1(S^1),\xi )$.
We then take this theorem as the starting point for the classification
of more general Legendrian `cable links' in $(J^1(S^1),\xi )$.
A further purpose of the present paper is to address a subtle
point in our earlier classification of Legendrian cable links
in~$S^3$, where previously we did not provide full details,
see~\cite[p.~154]{dige07}. The relevant argument relies in an
essential way on the classification of Legendrian
circular helix links.

Before stating our results, we need
to recall the definition of the classical invariants of
a Legendrian knot in $(J^1(S^1),\xi )$.
\section{Legendrian knots in $J^1(S^1)$}
\label{section:knots-jetspace}
By identifying $S^1$ with $\R /2\pi\Z$, we can visualise
a Legendrian knot $K\subset J^1(S^1)$ in its front projection to
a strip $[0,2\pi ]\times\R$ in the $xz$--plane. The
{\bf Thurston--Bennequin invariant} of $K$ is
\[ \tb (K)=\mbox{\rm writhe}(K)-\frac{1}{2}\# (\mbox{\rm cusps}(K)),\]
where the quantities on the right are computed from the front projection
of~$K$. This is the signed number of crossing changes required to
`unlink' $K$ from its push-off $K'$ in the $z$--direction (i.e.\
transverse to~$\xi$), that is, the number of crossing changes that
will allow one to separate the two knots in $J^1(S^1)$.

The {\bf rotation number} of the {\em oriented} Legendrian knot
$K$ is defined in terms of the front
projection as
\[ \rot (K)=\frac{1}{2}(c_- - c_+),\]
with $c_{\pm}$ the number of cusps oriented upwards or downwards,
respectively.

There are definitions of these invariants that do not rely on the
front projection, and which show that $\tb$ and $\rot$ are in fact
Legendrian isotopy invariants, cf.~\cite{dige07}.

The discussion in \cite{dige07} also shows that the Thurston--Bennequin
inequality for Legendrian realisations $L_0$ of the topological knot type
$\Lambda_0$ takes the form
\[ \tb (L_0)+|\rot (L_0)|\leq 0.\]
It follows that the maximal Thurston--Bennequin invariant of such knots $L_0$
is $\otb (L_0 )= 0$, and this maximal value is realised by the Legendrian
knot~$\Lambda_0$. By stabilising $\Lambda_0$, as described in the next
section, one can realise any non-positive integer $-m$ in the form
$-m=\tb (L_0)$.

In \cite[Section~6]{dige07}
we wrote down an explicit contact embedding $f\co J^1(S^1)\rightarrow S^3$
(both manifolds equipped with their standard contact structure). The
image of $f$ is the complement of a standard Legendrian unknot~$K_0$;
the knots $f(\Lambda_0)$ and $K_0$ form a Hopf link. We shall
frequently use the contactomorphism $f\co J^1(S^1)\rightarrow S^3
\setminus K_0$ to relate Legendrian links in the two manifolds.
\section{Legendrian helix links}
\label{section:helix}
\begin{figure}[h]
\centering
\includegraphics[scale=0.3]{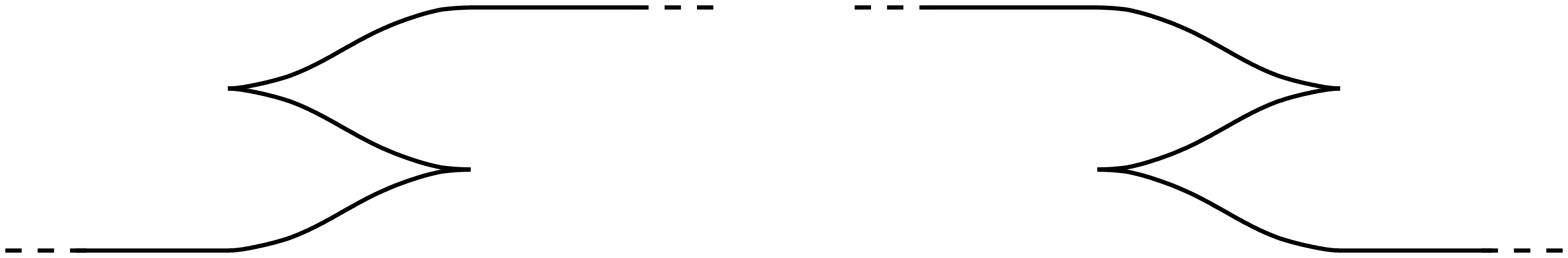}
  \caption{Legendrian `zigzags'.}
  \label{figure:zigzag}
\end{figure}

The {\bf stabilisation} $S_{\pm}(K)$ of a Legendrian knot $K$
in $(J^1(S^1),\xi )$ is the Legendrian knot whose front projection
is obtained from that of $K$ by adding a `zigzag' as in
Figure~\ref{figure:zigzag}.

A zigzag oriented downwards gives a positive stabilisation $S_+$,
while a zigzag oriented upwards gives~$S_-$. Hence
\[ \tb (S_{\pm}(K))=\tb (K)-1,\;\;\rot (S_{\pm}(K))=\rot (K)\pm 1.\]
Stabilisations are well defined and commute with each other.

\begin{thm}
\label{thm:helix}
Let $L=L_0\sqcup L_1$ be a Legendrian link in $(J^1(S^1),\xi )$
of the same (topological) oriented link type as $\Lambda_0\sqcup \Lambda_1$.

(a) If $\tb (L_0)=\tb (L_1)=0$, then $L$ is Legendrian isotopic
to either $\Lambda_0\sqcup \Lambda_1$ or $\Lambda_1\sqcup \Lambda_0$.

(b) If one of $\tb (L_0),\tb (L_1)$ is negative,
then $L$ is Legendrian isotopic to
\[ S_+^{k_0}S_-^{l_0}(\Lambda_0)\sqcup S_+^{k_1}S_-^{l_1}(\Lambda_1),\]
where
\[ k_i=\frac{-\tb (L_i)+\rot (L_i)}{2},\;\;
l_i=\frac{-\tb (L_i)-\rot (L_i)}{2},\;\; i=0,1.\]
\end{thm}

In particular, the phenomenon observed by Traynor disappears with
the first stabilisation. Theorem~\ref{thm:helix} will be proved by
a method analogous to that in \cite{etho01,dige07}.
In \cite{dige07}, one of the key ingredients was the
classification of tight contact structures on a thickened torus
$T^2\times [0,1]$, here it is the classification of tight contact
structures on the product of a pair of pants with a circle,
as obtained by Etnyre and Honda~\cite{etho01a}. The other main ingredient
is convex surface theory in the sense of Giroux; the exposition
given in \cite[Section~3]{etho01} is sufficient for our purposes.

\begin{proof}
(a) For $R>0$, consider the solid torus
\[ M_R=\{ (x,y,z)\in  J^1(S^1)\co y^2+z^2\leq R^2\} .\]
As meridian $\mu$ of this solid torus we may take the curve
on $\partial M_R$ given by $x=0$, oriented positively in
the $yz$--plane; the curve given by $(y,z)=(R,0)$, oriented by the
parameter~$x$, can serve as longitude~$\lambda$. Then $(\mu,\lambda )$
is a positive basis for $H_1(\partial M_R)$.

The vector field $X:=y\partial_y+z\partial_z$ satisfies
$\mathcal{L}_X(dz-y\, dx)=dz-y\, dx$, so it is a contact vector field
for~$\xi$. Hence $\partial M_R$ is a convex
torus whose dividing set --- the set of points where $X\in\xi$ ---
consists of two longitudes $\{ y=\pm R, z=0\}$;
i.e.\ these dividing curves have slope
$\infty$ with respect to $(\mu ,\lambda )$. The characteristic foliation
on $\partial M_R$ is the Legendrian ruling of slope~$1$ defined by
the vector field $\partial_x+y\partial_z-z\partial_y$; this foliation has
singularities along the Legendrian divides $\{ y=0, z=\pm R\}$.

Since $\tb (\Lambda_0)=\tb (\Lambda_1)=0$, we find disjoint tubular
neighbourhoods $N_0,N_1$ of these two knots with convex boundary
having two dividing curves of slope~$\infty$. Likewise, we may
choose disjoint tubular neighbourhoods $V_0,V_1$ of $L_0,L_1$,
respectively, with convex boundary of the same kind.
These neighbourhoods are contactomorphic by a diffeomorphism
preserving the longitude (and, obviously, the meridian). We may assume
that $R$ has been chosen large enough such that
these four tubular neighbourhoods are contained in $\Int (M_R)$. 

Write $\Sigma$ for a pair of pants, i.e.\ a disc with two open discs
removed. Both $M_R\setminus\Int (N_0\cup N_1)$ and
$M_R\setminus\Int (V_0\cup V_1)$ are homeomorphic to $\Sigma\times S^1$.
All the boundary tori are convex of slope $\infty$ with respect to $\xi$.
As shown in Lemma~10 and Lemma~11 of~\cite{etho01a}, this information
determines the tight contact structure on $\Sigma\times S^1$ up
to a permutation of the two interior `holes'; in other words,
we find a contactomorphism
\[ \phi\co M_R\setminus\Int (V_0\cup V_1)\longrightarrow
M_R\setminus \Int (N_0\cup N_1)\]
that is the identity on $\partial M_R$,
and sends $(\partial V_0,\partial V_1)$
to $(\partial N_{s_0},\partial N_{s_1})$ for some permutation
$(s_0,s_1)$ of $(0,1)$,
preserving the corresponding longitudes and meridians.

We can extend $\phi$ to a contactomorphism of $(J^1(S^1),\xi )$
that sends $L_0\sqcup L_1$ to $\Lambda_{s_0}\sqcup\Lambda_{s_1}$ and
equals the identity map outside~$M_R$. As explained
in the proof of \cite[Thm.~3.3]{dige07}, such a contactomorphism
is contact isotopic to the identity. This yields the Legendrian isotopy
of $L_0\sqcup L_1$ to $\Lambda_{s_0}\sqcup\Lambda_{s_1}$.

(b) Suppose that $\tb (L_0)<0$. Choose an annulus $A$, disjoint
from $L_1$, with boundary curves
$L_0$ and $\Lambda_s$ for some large constant~$s$. The surface framing
of $\Lambda_s$ determined by $A$ corresponds to the vector
field~$\partial_z$, so does the contact framing of $\Lambda_s$, which
satisfies $\tb (\Lambda_s)=0$. Since $\tb (L_0)<0$, the contact framing
of $L_0$ determined by $\xi$ makes at least one negative
(i.e.\ left) twist relative to its surface framing determined by~$A$.
The condition that the contact structure $\xi$
twist non-positively relative to the surface framing along either boundary
component of $A$ is precisely what is needed so that one may invoke
\cite[Prop.~3.1]{hond00}, which tells us that we can find a $C^0\!$--small
perturbation of $A$ rel $L_0,\Lambda_s$ (keeping it disjoint
from~$L_1$) into a convex surface.

In slightly less explicit form, this perturbation statement
is also part of the following result of Kanda~\cite{kand98}, cf.\
\cite[Thm.~4.8]{dige07}, to which we shall refer frequently in
the sequel.

\begin{thm}[Kanda]
\label{thm:Kanda}
If $\gamma$ is a Legendrian curve in a surface $\Sigma$, then $\Sigma$
may be isotoped relative to $\gamma$ so that it is convex if and only if
the twisting $t_{\Sigma}(\gamma )$ of the contact planes along $\gamma$
relative to the framing induced by $\Sigma$ satisfies $t_{\Sigma}(\gamma )
\leq 0$. If $\Sigma$ is convex, then
\[ t_{\Sigma}(\gamma )=-\frac{1}{2}\# (\gamma\cap\Gamma ),\]
where $\# (\gamma\cap\Gamma )$ denotes the number of intersection points
of $\gamma$ with the dividing set $\Gamma$ of~$\Sigma$.
\end{thm}

By this theorem,
the dividing curves of the (convex) annulus $A$ intersect $L_0,\Lambda_s$
in $-2\tb (L_0)>0$, $-2\tb (\Lambda_s)=0$ points, respectively.
This means that there is a boundary parallel dividing curve
along~$L_0$, which allows us to destabilise this Legendrian knot in
$J^1(S^1)\setminus L_1$, cf.\ \cite[Lemma~3.9]{etho01}.

By iterating this procedure we find Legendrian knots
$L_0', L_1'$ with $\tb (L_0')=\tb (L_1')=0$ and
non-negative integers $k_i,l_i$ such that $S_+^{k_i}S_-^{l_i}(L_i')
=L_i$, $i=0,1$. So, by (a), $L$ is Legendrian isotopic
to
\[ S_+^{k_0}S_-^{l_0}(\Lambda_{s_0})\sqcup S_+^{k_1}S_-^{l_1}(\Lambda_{s_1})\]
for some permutation $(s_0,s_1)$ of $(0,1)$.
From the behaviour of the classical invariants
$\tb$ and $\rot$ under stabilisation it follows that the
relation between $k_i,l_i$ and $\tb (L_i),\rot (L_i)$
is as stated in the theorem.

\begin{figure}[h]
\labellist
\small\hair 2pt
\pinlabel $0$ [B] at 182 -20
\pinlabel $2\pi$ [B] at 290 -20
\pinlabel $-\pi$ [B] at 361 -20
\pinlabel $\pi$ [B] at 469 -20
\endlabellist
\centering
\includegraphics[scale=0.45]{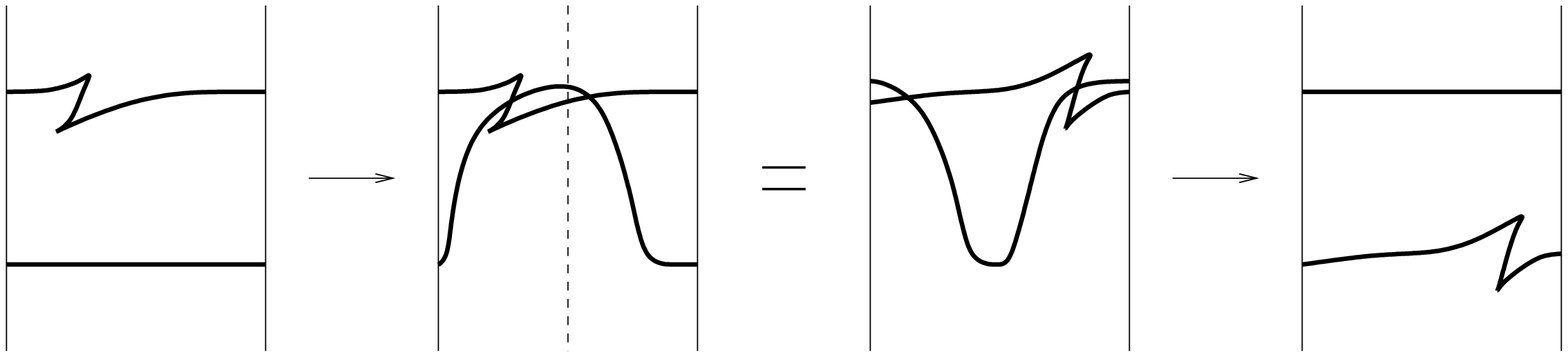}
\caption{A Legendrian isotopy between $\Lambda_0\sqcup S_+(\Lambda_1)$ and
$\Lambda_1\sqcup S_+(\Lambda_0)$.}
\label{figure:L1L2}
\end{figure}

It remains to show that in fact either permutation gives the
same Legendrian isotopy class, provided there is at least one
stabilisation. For that, it suffices to describe a Legendrian isotopy
between, say, $\Lambda_0\sqcup S_+(\Lambda_1)$ and
$\Lambda_1\sqcup S_+(\Lambda_0)$;
the picture for a single negative stabilisation is analogous.
Figure~\ref{figure:L1L2} indicates that a couple of Legendrian
Reidemeister moves of the second kind will do the job.
\end{proof}

\section{Cable links in $S^3$ revisited}
\label{section:S3}
In our earlier paper \cite{dige07} we gave a classification
of Legendrian $(p,q)$--cable links in~$S^3$. As pointed out
on p.~154 of that paper, we still owe the reader a final
detail in one particular case of that classification.
Here we provide the missing argument, which depends
crucially on Theorem~\ref{thm:helix}.

\begin{figure}[h]
\labellist
\small\hair 2pt
\pinlabel $\lambda$ at 342 36
\pinlabel $\mu$ at 208 103
\pinlabel $L_0$ at 218 63
\pinlabel $T$ [bl] at 460 210
\endlabellist
\centering
\includegraphics[scale=0.4]{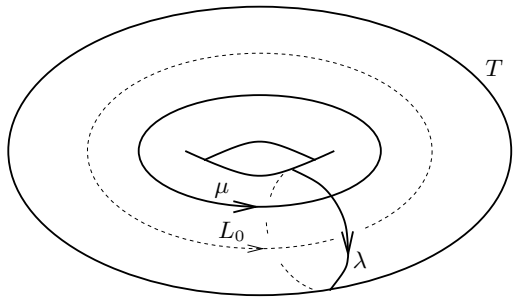}
\caption{Choice of meridian $\mu$ and longitude $\lambda$.}
\label{figure:torus}
\end{figure}

Recall the set-up from~\cite{dige07}, cf.\ Figure~\ref{figure:torus}.
Let $L_0$ be a trivial knot in $S^3$. Denote by $\mu$ and $\lambda$
meridian and longitude on a torus $T$, viewed as the boundary of
the {\em complement} of a tubular neighbourhood of~$L_0$; this unusual
convention was chosen in \cite{dige07} for technical reasons.
A {\bf $(p,q)$--cable link} is a link of the form $\LL=L_0\sqcup L_1$ with
$L_1$ a knot on $T$ homologically equivalent to $p\mu +q\lambda$.

In \cite{dige07} we classified Legendrian realisations of this
link type. With the notation $m:=-\tb (L_0)$,
the case for which additional details have to be provided is
the one where $q\geq 2$, $p=-1$ and $-m+q\leq 0$ (case 3(b2-iii)
in~\cite{dige07}). Here the maximal possible value of $\tb (L_1)$
is $pq$. What remains to be shown in order to complete the
argument from \cite{dige07} is the uniqueness of the Legendrian
realisation of this link type with $\tb (L_1)=pq$ and given rotation
numbers in the allowable range determined by equations (4.3) and
(4.4) of~\cite{dige07} (as we shall see, no explicit reference to
these equations will be necessary in the discussion that follows).

As explained in \cite[p.~147]{dige07}, we may assume
that the torus $T$ on which $L_1$ sits is a convex torus in standard
form (cf.~\cite[Defn.~4.1]{dige07}), and $L_1$ is one of an even number
of Legendrian divides on~$T$. The slope of these curves --- with $T$ now
regarded as the boundary of a tubular neighbourhood of~$L_0$ --- is
$p/q=-1/q$.

For such slopes, \cite[Lemma~3.17]{etho01} does not apply. For other slopes,
this lemma could be used to show that all Legendrian divides are isotopic.
In particular, one could use bypasses in order to remove all but two
Legendrian divides and still assume without loss of generality
that $L_1$ was one of these remaining Legendrian divides.

Instead, we argue as follows. The torus $T$ being convex, there
is a contact flow transverse to~$T$. Use this flow to push $T$
further away from~$L_0$, so that $\LL$ may now be regarded
as a link in the interior of a solid torus $V$ with convex boundary of
slope~$-1/q$. If there are more than two dividing curves (and hence Legendrian
divides) on $\partial V$ we can use bypasses on meridional discs
of the complementary solid torus for reducing the number of dividing curves
down to~$2$. Now ignoring the complement of~$V$, we can apply a
diffeomorphism to~$V$ (changing the longitude(!) of $\partial V$
from $-\mu$ to $-\mu +q\lambda$) such that the boundary slope becomes~$\infty$.

This allows us to identify $V$ with a standard neighbourhood of $\Lambda_0$
in $J^1(S^1)$. With respect to this identification, $\LL$ is
topologically isotopic to $\Lambda_0\sqcup\Lambda_1$. If $-m+q<0$
then we are done by Theorem~\ref{thm:helix}~(b), because in this situation
Legendrian realisations of this link are unique even inside~$V$.

If $-m+q=0$, then by Theorem~\ref{thm:helix}~(a) our link $\LL$
is Legendrian isotopic inside $V$ to one of $\Lambda_0\sqcup\Lambda_1$
or $\Lambda_1\sqcup\Lambda_0$. Observe that the two components of this
link are Legendrian push-offs of each other. Translated back into
$S^3$, this means that $\LL$ consists of a
topologically trivial Legendrian knot and
its Legendrian push-off, or what is called the {\bf $2$--copy} of
a topologically trivial
Legendrian knot in~\cite{mish03}. (For the orientations to be
consistent with the interpretation of $L_1$ as the push-off of~$L_0$,
however, we need to replace $L_1$ by its reverse~$\oL_1$.)
As a consistency check, observe that
$\tb (L_0)=-m=-q=pq=\tb (\oL_1)$
and
$\lk (L_0,\oL_1)=-q=\tb (L_0)$.
Moreover, formul\ae\ (4.3) and (4.4) in \cite{dige07} give
$\rot(L_0)=\rot(\oL_1)$.

As shown in \cite[Prop.~4.2a,b]{mish03}, cyclic permutations
(by a Legendrian isotopy) of the $N$--copy of a
topologically trivial Legendrian knot
are possible in~$S^3$ (in contrast with Theorem~\ref{thm:helix}~(a)).
This proves the uniqueness of the Legendrian realisations of $\LL$ in
this case and completes the argument from~\cite{dige07}.

\begin{rem}
According to \cite[Thm.~5.1a]{mish03}, only cyclic permutations of the
$N$--copy of the Legendrian unknot are possible. For $N=3$,
this result is equivalent to Theorem~\ref{thm:Traynor}. Indeed,
with $f\co J^1(S^1)\rightarrow S^3\setminus K_0$ being the contactomorphism
mentioned earlier, the Legendrian link $f(\Lambda_0)\sqcup f(\Lambda_1)
\sqcup K_0$ is the $3$--copy of the Legendrian unknot. Thus,
if one were able to permute $\Lambda_0$ and $\Lambda_1$ in $J^1(S^1)$,
one would have a transposition of $f(\Lambda_0)$ and
$f(\Lambda_1)$, even with $K_0$ fixed during the isotopy.
Conversely, if one could transpose $f(\Lambda_0)$ and $f(\Lambda_1)$
by a Legendrian isotopy $\phi_t$ that ends up moving $K_0$ to itself,
one could also find a contactomorphism of $S^3$ that transposes
$f(\Lambda_0)$ and $f(\Lambda_1)$, and fixes a neighbourhood of
$K_0$ pointwise. (Define the contactomorphism to be the time--1--map $\phi_1$
of the isotopy outside a standard tubular neighbourhood $N_0$ of~$K_0$, set
it equal to the identity on a smaller standard tubular neighbourhood
contained in the interior of both $N_0$ and $\phi_1(N_0)$,
and extend this to a contactomorphism on all of $S^3$ using the
uniqueness up to contactomorphism of the non-rotative tight contact
structure on $T^2\times [0,1]$ with two dividing curves on each
boundary component, as proved in \cite[Prop.~4.9]{hond00}.)
Our argument in \cite[Proof of Thm.~3.3]{dige07} then
shows that one could also find a Legendrian isotopy that
transposes $f(\Lambda_0)$ and $f(\Lambda_1)$ and keeps a neighbourhood of
$K_0$ fixed during the isotopy. Such an isotopy would induce a
permutation of $\Lambda_0$ and $\Lambda_1$ in~$J^1(S^1)$.
\end{rem}

\section{Cable links in $J^1(S^1)$}
Recall the definition of $M_R$ and $\mu,\lambda$ given at the
beginning of the proof of Theorem~\ref{thm:helix}.
Any torus of the form $T=\partial M_R$ for some $R>0$ is
called a {\bf standard torus}.
By a {\bf $(p,q)$--cable link} in $J^1(S^1)$ we mean a link
$\LL =L_0\sqcup L_1$ isotopic (as an ordered, oriented link)
to the union of $\Lambda_0$ and a $(p,q)$--cable of~$\Lambda_0$, i.e.\
a knot on a standard torus in the class of $p\mu +q\lambda$
(with $p$ and $q$ coprime). Our main theorem gives a
classification of the Legendrian realisations of this link type.

\begin{thm}
\label{thm:cable}
With the exception of the case in Theorem~\ref{thm:helix}~(a),
two oriented Legendrian cable links in $(J^1(S^1),\xi )$ are Legendrian
isotopic if and only if their oriented link types and their
classical invariants $\tb$ and $\rot$ agree.
\end{thm}

For a given pair $(p,q)$, let $\LL=L_0\sqcup L_1$ and
$\LL^*=L_0^*\sqcup L_1^*$
be two Legendrian realisations of the $(p,q)$--cable link having
the same classical invariants. According to \cite[Thm.~3.3]{dige07},
Legendrian torus knots in $(J^1(S^1),\xi )$ are determined by the
classical invariants. Applied to $L_0$ and $L_0^*$, regarded
as a $(0,1)$--torus knot, this result tells us
that we may assume $L_0=L_0^*=S^{k_0}_+S^{l_0}_-\Lambda_0$,
where $m:=k_0+l_0=-\tb (L_0)$ and $k_0-l_0=\rot (L_0)$. 

The strategy for proving Theorem~\ref{thm:cable} is now parallel
to~\cite{dige07}. First determine an upper bound $\otb (L_1)$ on
the Thurston--Bennequin invariant of
Legendrian realisations $L_1$ of $(p,q)$--cables of the given $L_0$, and
show that Legendrian realisations with non-maximal $\tb$ destabilise.
Show further that Legendrian realisations with maximal $\tb$ are determined
by $\rot$. If there are distinct Legendrian realisations with maximal $\tb$,
one also needs to understand the relationship
between their stabilisations. This last issue
will only arise in one subcase, where it can essentially be
settled by a reference to~\cite{dige07}.

By choosing the orientation of $L_1$
appropriately, we may assume $q\geq 0$ (and $p=1$ for $q=0$).

\vspace{2mm}

{\bf Case 1: \boldmath{$p=0$} and \boldmath{$q=1$}.}
This is the case dealt with in Theorem~\ref{thm:helix}.

\vspace{2mm}

{\bf Case 2: \boldmath{$p=1$} and \boldmath{$q=0$}.}
Here $L_1$ is a meridian of $L_0$ and
thus in particular a trivial knot. So the Thurston--Bennequin inequality
gives $\otb (L_1)=-1$, and this upper bound is obviously realised
by a standard Legendrian unknot linked once with~$L_0$.

By Theorem~\ref{thm:Kanda} we may assume that $L_1$ lies on
a convex {\bf standardly embedded} torus~$T$, i.e.\ a torus isotopic in
$J^1(S^1)\setminus L_0$ to a standard torus containing $L_0$ in
the interior. Choose $R>0$ large enough such that the standard torus
$T_{\infty}:=\partial M_R$ of slope $\infty$ contains $T$ in the
interior. Since $\tb (L_0)=-m$, we find a small tubular neighbourhood
of $L_0$, contained in the interior of~$T$, with convex boundary
having two dividing curves of slope~$-1/m$. Since slopes of
convex standardly embedded tori decrease as we
move away from $L_0$, and the slope can never be zero by the
tightness of~$\xi$, we see that the slope of $T$ must be of the
form $-r/s$ with coprime $r\in\N$, $s\in\N_0$, and $-r/s\leq -1/m$.
By Giroux flexibility, cf.\ \cite[Thm.~3.4]{hond00}, we may
perturb $T_{\infty}$
so that it still has two dividing curves of slope $\infty$,
but now a Legendrian ruling of slope~$0$.

Let $A$ be an annulus connecting $L_1$ and a Legendrian ruling curve
$L_{\infty}$ on~$T_{\infty}$. We have $t_A(L_1)=t_T(L_1)$ and
$t_A(L_{\infty})=t_{T_{\infty}}(L_{\infty})$, so by
Kanda's theorem --- and using the convexity of both $T$ and~$T_{\infty}$
--- we may assume that $A$ is convex.
If $2n$ denotes the number of dividing curves
on~$T$, the algebraic intersection number of $L_1$ with the dividing set
$\Gamma_T$
on $T$ equals $2nr$. Thus, if $n>1$ or $r>1$, there is a boundary
parallel dividing curve on $A$ along~$L_1$, which allows us to destabilise
that knot (obviously this can only happen if $\tb (L_1)$ was smaller than
$-1$ to start with). If $n=r=1$, but the geometric intersection
number between $L_1$ and $\Gamma_T$ is larger than the algebraic intersection
number, we can again destabilise. If $n=r=1$ and $\# (L_1\cap\Gamma_T)=2$,
then
$\tb (L_1)=t_A(L_1)=t_T(L_1)=-1$
by Theorem~\ref{thm:Kanda}, so we are in the case $\tb (L_1)=\otb (L_1)$.
The dividing set of $A$ then consists of two curves connecting $L_1$
with $L_{\infty}$, so by Giroux flexibility we may assume
that $A$ is foliated by Legendrian circles parallel to the boundary.
This Legendrian ruling defines a Legendrian isotopy between $L_1$
and $L_{\infty}$ (in the complement of~$L_0$), which proves uniqueness
of the Legendrian link~$\LL$.

\vspace{2mm}

{\bf Case 3: \boldmath{$p\geq 1$} and \boldmath{$q\geq 1$}.}
Given a Legendrian link in $J^1(S^1)$, we are going to study it via
its image in~$S^3\setminus K_0$ under the contactomorphism $f$
mentioned in Section~\ref{section:knots-jetspace}.
We label objects in $S^3$ with a prime, and objects in $J^1(S^1)$ without.
It will be convenient (for comparison with the
classification in~\cite{dige07}) to define $L_0':=f(L_0)$ and
$L_1':=\overline{f(L_1)}$ (i.e.\ $f(L_1)$ with reversed
orientation). It is a straightforward check that if $\LL =L_0\sqcup L_1$
is a $(p,q)$--cable link
in $J^1(S^1)$, then $\LL':=L_0'\sqcup L_1'$ is a $(p',q')$--cable link
in $S^3$, where $p'=-q$ and $q'=q-p$. (Here meridian $\mu'$ and longitude
$\lambda'$ on a torus $T'$ around $L_0'$ are chosen as in
Figure~\ref{figure:torus}.) With $m:=-\tb (L_0)$ and
$m':=-\tb (L_0')$ we have $-m=-m'+1^2=-m'+1$.
(See \cite[p.~138]{dige07} and the proof of Proposition~\ref{prop:max-tb}
below for the relation between the Thurston--Bennequin invariant
in $J^1(S^1)$ and that same invariant in~$S^3$.)

However, we shall see that the classification
of Legendrian cable links in $S^3$ does not translate directly into a
classification of Legendrian cable links in~$J^1(S^1)$.

\begin{prop}
\label{prop:max-tb}
If $p\geq 1$ and $q\geq 1$, then $\otb (L_1)=p(q-1)$.
\end{prop}

\begin{proof}
In Figure~\ref{figure:p-q-link1} we exhibit the front projection of
a $(p,q)$--cable $L_1$ to $\Lambda_0$ (or any of its stabilisations)
with $p,q\geq 1$ and $\tb (L_1) = p(q-1)$
(which can be read off as the writhe of that front).
This example shows that $\otb (L_1)\geq p(q-1)$.

\begin{figure}[h]
\labellist
\small\hair 2pt
\pinlabel {$q$ strands} [r] at 0 144
\pinlabel {$p$ such crossings} at 90 10
\endlabellist
\centering
\includegraphics[scale=0.4]{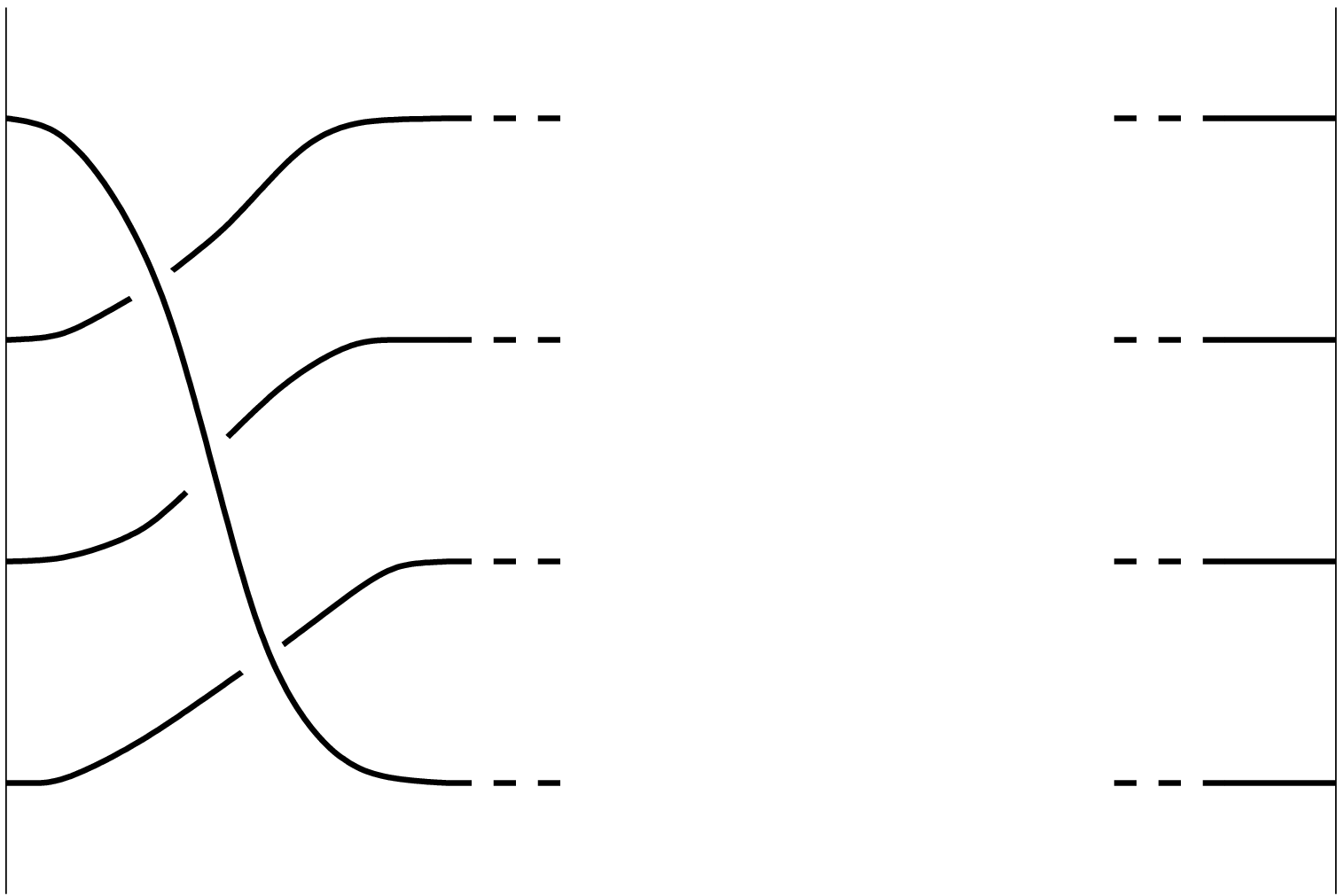}
\caption{$(p,q)$--cable with maximal $\tb$ in Case 3.}
\label{figure:p-q-link1}
\end{figure}

It remains to prove the complementary inequality $\tb (L_1)\leq p(q-1)$
for every Legendrian $(p,q)$--cable $L_1$ of~$L_0$.

If $p\geq q$, then the $(-p',-q')$--torus knot $f(L_1)$
is a {\em positive} torus knot.
The maximal Thurston--Bennequin invariant of a positive $(-p',-q')$--torus knot
is, according to \cite[Thm.~4.1]{etho01}, given by $p'q'+p'+q'$
(this also holds for $q'=0$, in which case $p'=-1$).
With the correction term $q^2$ --- for a knot homotopic to $q$ times
the generator of $\pi_1(J^1(S^1))$ --- between the Thurston--Bennequin
invariant in
$J^1(S^1)$ and that in $S^3$ (cf.\ \cite[p.~138]{dige07}), we find
\[ \tb (L_1) \leq (p-q)q-(p-q)-q+q^2=p(q-1).\]

If $0<p<q$, we consider the
$(-q',-p')$--cable link $K_0\sqcup f(L_1)$; here $(-q',-p')=(p-q,q)$
with $q\geq 2$. Thus, we produce
a cable link in $S^3$ corresponding to Case~3(b) of~\cite{dige07}.
Then \cite[Lemma~4.6]{dige07} tells us that $\tb (f(L_1))\leq p'q'$.
(Observe that the situation discussed there, where we might have
to interchange the roles of the two link components, does not
arise: this is a consequence of $\tb (K_0)=-1>-q$; see the discussion
around \cite[Lemma~4.5]{dige07}.)
We compute
\[ \tb (L_1)=\tb(f(L_1))+q^2\leq -q(q-p)+q^2=pq.\]

So in all cases $p,q\geq 1$ we have at least the estimate $\tb (L_1)
\leq pq$, and the following argument applies.
(As regards the proof of the proposition, this argument is of course
superfluous if $p\geq q$, but later it will be employed to conclude the
classification in Case~3, where we need it for all $p,q\geq 1$.)
Let $T$ be a standardly embedded torus on which $L_1$ sits.
Then $t_T(L_1)=\tb (L_1)-pq\leq 0$, so Theorem~\ref{thm:Kanda}
allows us to assume that $T$ is convex. As in Case~2, the dividing set
$\Gamma_T$ consists of $2n$ curves of slope $-r/s\leq -1/m$. Thus
\[ \# (L_1\cap\Gamma_T)\geq 2n\left|\begin{array}{rr}p&-s\\q&r
\end{array}\right|\geq 2p.\]
Hence $t_T(L_1)\leq -p$ (again by Theorem~\ref{thm:Kanda}), and therefore
$\tb (L_1)=pq+t_T(L_1)\leq p(q-1)$, as claimed.
\end{proof}

\begin{rem}
For $0<p<q$, the estimate $\tb (L_1)\leq p(q-1)$ can also be obtained by
appealing directly to \cite[Lemma~4.6]{dige07}, applied to
the $(p-q,q)$--cable link $K_0\sqcup f(L_1)$. Since $\tb (K_0)=-1$,
that lemma yields
\[ \tb (L_1)=\tb (f(L_1))+q^2\leq (p-q)q-\max (1\cdot(p-q)+q,0)+q^2=p(q-1).\]
\end{rem}

Now the reasoning for Case~3 is entirely analogous to Case~2.
We take $T_{\infty}$
to be a (sufficiently large) convex standardly embedded torus having two
dividing curves of slope~$\infty$,
and $L_{\infty}\subset T_{\infty}$ a ruling curve in a Legendrian ruling
of slope~$q/p$. Then $\# (L_{\infty}\cap\Gamma_{T_{\infty}})=2p$.

On the other hand, we have
$\# (L_1\cap\Gamma_T)\geq 2p$, with equality if and only if
$r=1$, $s=0$, $n=1$, and geometric intersection number being
equal to the algebraic one. In all other cases, we can destabilise~$L_1$.
If $\tb(L_1)=\otb (L_1)$, we must have $\# (L_1\cap\Gamma_T)=2p$,
and with the help of a convex
annulus $A$ between $L_1$ and $L_{\infty}$ we show those two curves
to be Legendrian isotopic as in Case~2. This concludes the
proof of Theorem~\ref{thm:cable} in Case~3.

\vspace{2mm}

The next corollary shows that in some cases a
Legendrian $(p',q')$--cable link in $S^3\setminus K_0$ can
be further destabilised in~$S^3$, but not in $S^3\setminus K_0$.

\begin{cor}
\label{cor:noimage}
For $0<p<q-1$ there are Legendrian $(p',q')$--cable links in
$S^3$ that are not Legendrian isotopic to
the image under $f$ of a Legendrian $(p,q)$--cable link
in~$J^1(S^1)$.
\end{cor}

\begin{proof}
Let $L_0'\sqcup L_1'$ be a Legendrian $(p',q')$--cable link in $S^3$
with $q'\geq 2$, $p'<0$, and the unknotted component $L_0'$
having $\tb (L_0')=-m'$. According to \cite[Lemma~4.6]{dige07},
the maximal Thurston--Bennequin invariant $\otb (L_1')$ is given
by $p'q'-\max (m'p'+q',0)$. If $L_0'\sqcup L_1'=f (L_0)\sqcup
\overline{f(L_1)}$ is the
image of a Legendrian $(p,q)$--cable link in $J^1(S^1)$
with $\tb (L_0)=-m$, then
\[ m'p'+q'= (m+1)\cdot (-q)+q-p=-mq-p<0.\]
If we consider arbitrary Legendrian realisations $L_0^*\sqcup L_1^*$
of this link type in $S^3$ (with $\tb (L_0^*)=-m'$), not necessarily
coming from a link in $J^1(S^1)$, then
$\otb (L_1^*)=p'q' =pq-q^2$. This is larger than
\[ \otb (f(L_1))=\otb (L_1)-q^2=p(q-1)-q^2, \]
so a sufficiently high destabilisation of $L_1'$ in
$S^3\setminus L_0'$ gives a
Legendrian $(p',q')$--cable link which is not Legendrian isotopic
to the image of a Legendrian $(p,q)$--cable link in~$J^1(S^1)$.
\end{proof}

Thus, if $L_0\sqcup L_1$ is a
Legendrian $(p,q)$--cable link in $J^1(S^1)$ with
$0<p<q-1$ and $\tb (L_1)=\otb (L_1)=p(q-1)$, then
$\overline{f(L_1)}$ can be destabilised in $S^3\setminus f(L_0)$,
but not in $S^3\setminus (f(L_0)\cup K_0)$. Still, our uniqueness
proof above implies that there can be no choice for $\rot (L_1)$.
This may seem surprising, since the classification from~\cite{dige07},
applied to $f(L_0)\sqcup \overline{f(L_1)}$ in~$S^3$, allows
freedom in $\rot (\overline{f(L_1)})$ if $\tb (f(L_1))<\otb (L_1')$.
The paradox is resolved by considering, again, $K_0\sqcup f(L_1)$.
Then \cite[Lemma~4.11]{dige07} gives
\[ \rot (L_1)=\rot (f(L_1))=(p-q)\cdot\rot (K_0)=0.\]
This is indeed the rotation number realised by the cable in
Figure~\ref{figure:p-q-link1}.

\vspace{2mm}

{\bf Case 4: \boldmath{$p<0$} and \boldmath{$q\geq 1$}.}
Here we have $p'=-q<0$ and $q'=q-p\geq 2$, so we are in Case~3(b)
of~\cite{dige07}. Recall from the proof of Corollary~\ref{cor:noimage}
that $m'p'+q'=-mq-p$.

\vspace{2mm}

{\bf Case 4(a): \boldmath{$p<0$}, \boldmath{$q\geq 1$} and
\boldmath{$mq+p<0$}.}
This puts us in Case~3(b1) of~\cite{dige07}.
The Legendrian $(p',q')$--cable link $L_0'\sqcup L_1'$ in $S^3$
realising $\otb (L_1')$, constructed explicitly in Section~5 of our earlier
paper, sits in $S^3\setminus K_0$ and thus gives a Legendrian realisation
of a $(p,q)$--cable link in $J^1(S^1)$ with
\[ \tb (L_1)=\otb (L_1)= p'q'-(m'p'+q')+q^2=pq+mq+p.\]
Here is a more direct description of such a Legendrian cable link realising
$\otb(L_1)$. Consider a small standard neighbourhood of $L_0$
with convex boundary $T_0$ having two dividing curves
of slope $-1/m$. Since $q/p\neq -1/m$, Giroux flexibility allows us
to assume that $T_0$ has a Legendrian ruling of slope $q/p$.
Let $L_1^0$ be one of these ruling curves. Then, by Kanda's theorem,
\[ \tb (L_1^0)=pq+t_{T_0}(L_1^0)=pq-\left|\begin{array}{cr}p&m\\
q&-1\end{array}\right|=pq+mq+p=\otb(L_1).\]

The following argument is analogous to Cases 2 and~3, except that
the protagonist is now $T_0$ in place of~$T_{\infty}$.
Thus, let $L_1$ be any other Legendrian $(p,q)$--cable of~$L_0$,
lying on some standardly embedded torus~$T$. The neighbourhood
in the previous construction can be chosen arbitrarily small,
so we may assume that $T_0$ sits in the interior of~$T$. The fact that
\[ t_T (L_1)=\tb (L_1)-pq\leq\otb (L_1)-pq=mq+p<0\]
allows us to assume that $T$ is convex, having $2n$ dividing
curves of slope $-r/s\leq -1/m$.

Observe that, since $q>0$, we have $sq+rp\leq r(mq+p)<0$. Hence
\[ \# (L_1\cap\Gamma_T)\geq -2\left|\begin{array}{cr}p&-s\\
q&r\end{array}\right|\geq -2(mq+p)=\# (L_1^0\cap\Gamma_{T_0}),\]
with equality if and only if $\Gamma_T$ consists of two dividing curves of
slope $-r/s=-1/m$. With the help of a convex annulus between
$L_1^0$ and $L_1$ we conclude that $L_1$ can be destabilised for
$\tb (L_1)<\otb (L_1)$, and that $L_1$ is Legendrian isotopic to
$L_1^0$ if its $\tb$ is maximal.

\vspace{2mm}

{\bf Case 4(b): \boldmath{$p<0$}, \boldmath{$q\geq 1$} and
\boldmath{$mq+p\geq 0$}.}
This corresponds to Case~3(b2) of~\cite{dige07}. If $q=1$, we
are in subcase 3(b2-iii); if $q>1$, in subcase 3(b2-ii).
Here we have $\otb (L_1')=p'q'=pq-q^2$, hence $\tb(L_1)\leq pq$.
There are several Legendrian realisations of this maximal $\tb$,
distinguished by $\rot (L_1')$. All these realisations,
as described in \cite[Section~5]{dige07}, live in $S^3\setminus K_0$
in such a way that we may regard them as links in~$J^1(S^1)$.
This yields $\otb (L_1)=pq$. In order to establish
a one-to-one correspondence between Legendrian
$(p,q)$--cable links in $J^1(S^1)$ and $(p',q')$--cable links in~$S^3$,
we need to show the following:

\begin{itemize}
\item In a $(p',q')$--cable link in $S^3\setminus K_0$ with non-maximal
$\tb (L_1')$, the component $L_1'$
can be destabilised inside $S^3\setminus K_0$.
\item Two $(p',q')$--cable links in $S^3\setminus K_0$ with maximal
$\tb (L_1')$ and the same classical invariants are Legendrian
isotopic inside $S^3\setminus K_0$.
\item The relation between stabilisations with the same classical
invariants coming from links with $\tb(L_1')=\otb (L_1')$ but
different rotation numbers are as described in Lemmata 4.20 and 4.21
of~\cite{dige07}, respectively.
\end{itemize}

The first point we can settle directly in $J^1(S^1)$. If
$L_0\sqcup L_1$ is a $(p,q)$--cable link with $\tb (L_1)<pq$,
we find a convex standardly embedded torus $T$ on which $L_1$ sits.
If the actual geometric intersection number between $L_1$ and the dividing
set $\Gamma$ of $T$ is larger than the minimal one possible,
there exists a bypass for $L_1$ on $T$, and we can destabilise~$L_1$.
Thus, we may assume that $T$ is in standard form of slope $-r/s\leq
-1/m$, and the non-maximality of $\tb (L_1)$ then forces
$-r/s\neq q/p$. Further, the condition $mq+p\geq 0$ gives $q/p\leq -1/m$.
So we find a convex torus $T_{q/p}$ of slope $q/p$ either between $T_0$
and~$T$ (if $q/p>-r/s$) or between $T$ and~$T_{\infty}$ (if $q/p<-r/s$),
and we can use a convex annulus between $L_1$ and a Legendrian divide on
$T_{q/p}$ to destabilise~$L_1$.

The second point, for $q>1$, requires the use of \cite[Lemma~4.14]{dige07}.
The proof of that lemma goes through in $S^3\setminus K_0$,
as can be seen by thinking of $S^3$ as being composed of a
union of four pieces: a standard neighbourhood of $L_0$ with boundary~$T_0$,
a thickened torus between $T_0$ and~$T$ (on which $L_1$ sits),
a thickened torus between $T$ and $T_{\infty}$, and a standard neighbourhood
of~$K_0$.

For $q=1$, this point cannot be settled with a simple reference
to~\cite{dige07}, for this is exactly the situation that necessitated
the discussion in Section~\ref{section:S3}. That argument allows us to
assume that our cable link sits in the interior of a solid torus
$V$ with convex boundary having two dividing curves of slope $1/p$.
(The reduction of the number of dividing curves on $\partial V$
to $2$ by using meridional discs in the complement also works in
$S^3\setminus K_0=J^1(S^1)$, because we can always
choose the bypasses in such a way that we do not move across
$K_0$ in order to effect the destabilisation.)

By `unwinding' this solid torus $V$ as in Section~\ref{section:S3},
we may identify $L_0\sqcup L_1$ with a link 
$L_0^*\sqcup L_1^*$ in
$J^1(S^1)$ topologically isotopic to $\Lambda_0\sqcup\Lambda_1$,
with $\tb (L_0^*)=-m-p$ and $\tb (L_1^*)=p\cdot 1-p=0$. Thus,
if $m+p>0$, then this link is determined by its classical
invariants even inside~$V$. If $m+p=0$, then $L_0^*\sqcup L_1^*$ is
Legendrian isotopic to $\Lambda_0\sqcup\Lambda_1$ or
$\Lambda_1\sqcup\Lambda_0$, and in particular a $2$--copy. So $L_0\sqcup L_1$
is likewise a $2$--copy, but now of a stabilised knot. In that case,
it is easy to see that permutation is possible by a Legendrian isotopy,
see~\cite[Prop.~4.2a]{mish03}.

Finally, for the third point, a close inspection of the relevant arguments
in \cite{dige07} reveals that the mentioned lemmata do indeed hold
in $S^3\setminus K_0$.

This completes the discussion of Case~4 and hence the proof of
Theorem~\ref{thm:cable}.

{\footnotesize
\begin{ack}
Most of this research was carried out during a stay of F.~D.\ at the
Mathematical Institute of the Universit\"at zu K\"oln, supported by DAAD
grant A/06/27941. F.~D.\ is partially supported by grant no.\
10631060 of the National Natural Science Foundation of China.
H.~G.\ would like to thank the American Institute of Mathematics
for its support during the workshop ``Legendrian and transverse knots''
(September 2008), where the final writing was done, and the organisers
D.~Fuchs, S.~Tabachnikov and L.~Traynor for the invitation to that
workshop.
\end{ack}
}

\end{document}